\documentclass[11pt]{article}%
\usepackage{amsmath}
\usepackage{amsfonts}
\usepackage{amssymb}
\usepackage{graphicx}
\usepackage{theorem}
\usepackage{makeidx}
\usepackage{verbatim}%
\providecommand{\U}[1]{\protect\rule{.1in}{.1in}}
\newtheorem{thm}{Theorem}[section]
\newtheorem{lm}[thm]{Lemma}
\newtheorem{pr}[thm]{Proposition}
\newtheorem{df}[thm]{Definition}
\newtheorem{rmk}[thm]{Remark}
\newtheorem{cor}[thm]{Corollary}
{\theorembodyfont{\upshape}

}
\numberwithin{equation}{section} 
\begin{document}

\title{ }

\begin{center}
\vspace*{1.5cm}

\textbf{FILLING THE GAP BETWEEN METRIC REGULARITY AND FIXED POINTS: THE
LINEAR\ OPENNESS OF\ COMPOSITIONS}

\vspace*{1cm}

M. DUREA

{\small {Faculty of Mathematics, "Al. I. Cuza" University,} }

{\small {Bd. Carol I, nr. 11, 700506 -- Ia\c{s}i, Romania,} }

{\small {e-mail: \texttt{durea@uaic.ro}}}

\bigskip

R. STRUGARIU

{\small Department of Mathematics, "Gh. Asachi" Technical University, }

{\small {Bd. Carol I, nr. 11, 700506 -- Ia\c{s}i, Romania,} }

{\small {e-mail: \texttt{rstrugariu@tuiasi.ro}}}
\end{center}

\bigskip

\bigskip

\noindent{\small {\textbf{Abstract: }}This paper is devoted to the
investigation of an important issue recently brought into attention by a
recent paper of Arutyunov: the relation between openness of composition of
set-valued maps and fixed point results. More precisely, we prove a general
result concerning the openness of compositions and then we show that this
result covers and implies most of the known openness results. In particular,
we reobtain several recent results in this field, including a fixed point
theorem of Dontchev and Frankowska.}

\bigskip

\noindent{\small {\textbf{Keywords: }}composition of set-valued mappings
$\cdot$ linear openness $\cdot$ metric regularity $\cdot$ Lipschitz-like
property $\cdot$ implicit multifunctions $\cdot$ fixed points}

\bigskip

\noindent{\small {\textbf{Mathematics Subject Classification (2010): }90C30
$\cdot$ 49J53} {$\cdot$ 54C60}}

\bigskip

\section{Introduction}

The equivalent properties of metric regularity and openness at linear rate,
which are nowadays mainly studied in the context of the set-valued maps, have
an important history. Their origins are in the open mapping principle for
linear operators obtained in the 1930s by Banach and Schauder. Subsequently, a
nonlinear extension of open mapping principle was obtained by Lyusternik
(1934) and Graves (1950) and, besides the importance on this result, another
important contribution is its proof which became, along the years, a powerful
instrument in the effort of getting generalizations, extensions and a better
understanding of these seminal works. Another landmark was the extension of
research to the case of set-valued maps with closed and convex graph and was
done by Ursescu in 1975 and Robinson in 1976, respectively. The famous
Robinson-Ursescu Theorem and a series of works of Milyutin in 1970s concerning
the preservation of regularity (and openness) of set-valued maps under
functional perturbation give a strong impetus to this direction of research
and, therefore, since 1980s to our time many mathematicians have participated
into a joint effort on the development and understanding of these problems. We
mention here only a few, having major contributions to the field: J. P. Aubin,
A. Dontchev, H. Frankowska, A. Ioffe, B. S. Mordukhovich, J. P. Penot, R. T.
Rockafellar, S. M. Robinson, C. Ursescu. We remark that in several works (of
J. P. Penot and C. Ursescu, for instance) dealing with openness results, the
techniques based on the original proof of Lyusternik-Graves Theorem (which is
in fact an iteration procedure) have been replaced by the use of Ekeland
Variational Principle. Detailed accounts on the historical facts, as well as
on the evolution of terminology are given in the monographs of Rockafellar and
Wets \cite{RocWet}, Mordukhovich \cite{Mor2006} and Dontchev and Rockafellar
\cite{DontRock2009b}.

The present paper is motivated by some (very) recent developments, which
suggest that the study of metric regularity can be done in relation with some
fixed points results. This new approach is suggested in the monograph of
Dontchev and Rockafellar \cite{DontRock2009b}, and it is more precisely
emphasized in the work of Arutyunov \cite{Arut2007}. Subsequently, the
relation between fixed points results and metric regularity are further
investigated in \cite{DonFra2010} and \cite{Ioffe2010b}. In \cite{DonFra2010},
\cite{Ioffe2010b}, some compositions of set-valued maps are involved and, in
fact, the degrees of generality of considered compositions represent, in some
extent, the main novelties in both these papers.

More precisely, Ioffe \cite{Ioffe2010b} investigates separately coincidence
fixed point results and openness of compositions, while Dontchev and
Frankowska in \cite{DonFra2010} present a fixed point theorem which
generalized the previous result of Arutyunov, and, moreover, show that their
result implies several well-known results concerning the metric regularity of
sum of set-valued maps. Moreover, these papers emphasize the fact that the
link between fixed points and metric regularity of the involved multifunctions
has a purely metric behavior. As a consequence, their proofs make use of
well-fitted adaptations of the initial iterative procedure of
Lyusternik-Graves type, and work on general metric spaces, under local
completeness assumptions.

In this paper we propose ourselves to bring into light a global approach and
to fill several implications which are not previously given. Having this aim
in mind, we follow the next procedure. Firstly, we prove some implicit
multifunctions assertions which, besides their own importance are used several
times as auxiliary results. Secondly, we prove a general openness result for
compositions of set-valued maps. Our proof is a refinement of the method based
on Ekeland Variational Principle. In particular, we prove that from our result
one can get a nontrivial extension of the corresponding result of Ioffe.
Moreover, the same result implies as well Dontchev-Frankowska fixed point
theorem and this implication (i.e. from openness assertions to fixed points
results) is not given elsewhere (up to our knowledge). At this point, it is
clear that all the implications from Dontchev-Frankowska fixed point theorem
apply, and therefore our openness result of compositions implies several very
important results in this topic. At every step of our investigation we are
trying to give an accurate account on the usual difficulties that arise at
certain points of the proofs and to interrelate our results and arguments to
the ones in literature. The setting we assume in this work is that of Banach
spaces. Note that most of the results work equally well on normed vector
spaces, with corresponding completeness assumptions around the reference points.

The paper is organized as follows. In the second section we introduce the
basic notations, concepts and we recall some known results. The third section
contains the main results of the paper. We present here some implicit
multifunction theorems which are important as separate results, but their
conclusions are in force in the proof of other results. Next, we formulate and
prove the main result of the paper which is an openness theorem for
compositions of set-valued maps. Then, on this basis, one result of Ioffe is
significantly extended and completed. The last section is devoted to the proof
of Dontchev-Frankowska fixed point theorem as a consequence of our main result
and several interpretations and possibilities of generalizations are
presented. The paper ends with a short conclusion section.

\section{Preliminaries}

This section contains some basic definitions and results used in the sequel.
In what follows, we suppose that all the involved spaces are Banach. In this
setting, $B(x,r)$ and $D(x,r)$ denote the open and the closed ball with center
$x$ and radius $r,$ respectively$.$ Sometimes we write $\mathbb{D}_{X}$ for
the closed unit ball of $X$. If $x\in X$ and $A\subset X,$ one defines the
distance from $x$ to $A$ as $d(x,A):=\inf\{\left\Vert x-a\right\Vert \mid a\in
A\}.$ As usual, we use the convention $d(x,\emptyset)=\infty.$ For a non-empty
set $A\subset X$ we put $\operatorname*{cl}A$ for its topological closure.
When we work on a product space, we consider the sum norm, unless otherwise stated.

Consider now a multifunction $F:X\rightrightarrows Y$. The domain and the
graph of $F$ are denoted respectively by
\[
\operatorname*{Dom}F:=\{x\in X\mid F(x)\neq\emptyset\}
\]
and%
\[
\operatorname*{Gr}F=\{(x,y)\in X\times Y\mid y\in F(x)\}.
\]
If $A\subset X$ then $F(A):=%
{\displaystyle\bigcup\limits_{x\in A}}
F(x).$ The inverse set-valued map of $F$ is $F^{-1}:Y\rightrightarrows X$
given by $F^{-1}(y)=\{x\in X\mid y\in F(x)\}$.

Recall that a multifunction $F$ is inner semicontinuous at $(x,y)\in
\operatorname*{Gr}F$ if for every open set $D\subset Y$ with $y\in D,$ there
exists a neighborhood $U\in\mathcal{V}(x)$ such that for every $x^{\prime}\in
U,$ $F(x^{\prime})\cap D\neq\emptyset$ (where $\mathcal{V}(x)$ stands for the
system of the neighborhoods of $x$).

We remind now the concepts of openness at linear rate, metric regularity and
Lipschitz-likeness of a multifunction around the reference point.

\begin{df}
\label{around}Let $L>0,$ $F:X\rightrightarrows Y$ be a multifunction and
$(\overline{x},\overline{y})\in\operatorname{Gr}F.$

(i) $F$ is said to be open at linear rate $L>0,$ or $L-$open around
$(\overline{x},\overline{y})$ if there exist a positive number $\varepsilon>0$
and two neighborhoods $U\in\mathcal{V}(\overline{x}),$ $V\in\mathcal{V}%
(\overline{y})$ such that, for every $\rho\in(0,\varepsilon)$ and every
$(x,y)\in\operatorname*{Gr}F\cap\lbrack U\times V],$%
\begin{equation}
B(y,\rho L)\subset F(B(x,\rho)). \label{Lopen}%
\end{equation}

The supremum of $L>0$ over all the combinations $(L,U,V,\varepsilon)$ for
which (\ref{Lopen}) holds is denoted by $\operatorname*{lop}F(\overline
{x},\overline{y})$ and is called the exact linear openness bound, or the exact
covering bound of $F$ around $(\overline{x},\overline{y}).$

(ii) $F$ is said to be Lipschitz-like, or has Aubin property around
$(\overline{x},\overline{y})$ with constant $L>0$ if there exist two
neighborhoods $U\in\mathcal{V}(\overline{x}),$ $V\in\mathcal{V}(\overline{y})$
such that, for every $x,u\in U,$%
\begin{equation}
F(x)\cap V\subset F(u)+L\left\Vert x-u\right\Vert \mathbb{D}_{Y}.
\label{LLip_like}%
\end{equation}

The infimum of $L>0$ over all the combinations $(L,U,V)$ for which
(\ref{LLip_like}) holds is denoted by $\operatorname*{lip}F(\overline
{x},\overline{y})$ and is called the exact Lipschitz bound of $F$ around
$(\overline{x},\overline{y}).$

(iii) $F$ is said to be metrically regular around $(\overline{x},\overline
{y})$ with constant $L>0$ if there exist two neighborhoods $U\in
\mathcal{V}(\overline{x}),$ $V\in\mathcal{V}(\overline{y})$ such that, for
every $(x,y)\in U\times V,$%
\begin{equation}
d(x,F^{-1}(y))\leq Ld(y,F(x)). \label{Lmreg}%
\end{equation}

The infimum of $L>0$ over all the combinations $(L,U,V)$ for which
(\ref{Lmreg}) holds is denoted by $\operatorname*{reg}F(\overline{x}%
,\overline{y})$ and is called the exact regularity bound of $F$ around
$(\overline{x},\overline{y}).$
\end{df}

The next proposition contains the well-known links between the notions
presented above. For more details about the proof, see \cite[Theorems 1.49,
1.52]{Mor2006}.

\begin{pr}
\label{link_around}Let $F:X\rightrightarrows Y$ be a multifunction and
$(\overline{x},\overline{y})\in\operatorname{Gr}F.$ Then $F$ is open at linear
rate around $(\overline{x},\overline{y})$ iff $F^{-1}$ is Lipschitz-like
around $(\overline{y},\overline{x})$ iff $F$ is metrically regular around
$(\overline{x},\overline{y})$. Moreover, in every of the previous situations,%
\[
(\operatorname*{lop}F(\overline{x},\overline{y}))^{-1}=\operatorname*{lip}%
F^{-1}(\overline{y},\overline{x})=\operatorname*{reg}F(\overline{x}%
,\overline{y}).
\]

\end{pr}

It is well known that the corresponding "at point" properties are
significantly different from the "around point" ones. Let us introduce now
some of these notions. For more related concepts we refer to
\cite{ArtMord2009}.

\begin{df}
\label{at}Let $L>0,$ $F:X\rightrightarrows Y$ be a multifunction and
$(\overline{x},\overline{y})\in\operatorname{Gr}F.$

(i) $F$ is said to be open at linear rate $L,$ or $L-$open at $(\overline
{x},\overline{y})$ if there exists a positive number $\varepsilon>0$ such
that, for every $\rho\in(0,\varepsilon),$%
\begin{equation}
B(\overline{y},\rho L)\subset F(B(\overline{x},\rho)). \label{Lopen_at}%
\end{equation}

The supremum of $L>0$ over all the combinations $(L,\varepsilon)$ for which
(\ref{Lopen_at}) holds is denoted by $\operatorname*{plop}F(\overline
{x},\overline{y})$ and is called the exact punctual linear openness bound of
$F$ at $(\overline{x},\overline{y}).$

(ii) $F$ is said to be pseudocalm with constant $L,$ or $L-$pseudocalm at
$(\overline{x},\overline{y}),$ if there exists a neighborhood $U\in
\mathcal{V}(\overline{x})$ such that, for every $x\in U,$%
\begin{equation}
d(\overline{y},F(x))\leq L\left\Vert x-\overline{x}\right\Vert .
\label{Lpseudocalm}%
\end{equation}

The infimum of $L>0$ over all the combinations $(L,U)$ for which
(\ref{Lpseudocalm}) holds is denoted by $\operatorname*{psdclm}F(\overline
{x},\overline{y})$ and is called the exact bound of pseudocalmness for $F$ at
$(\overline{x},\overline{y}).$

(iii) $F$ is said to be metrically hemiregular with constant $L,$ or
$L-$metrically hemiregular at $(\overline{x},\overline{y})$ if there exists a
neighborhood $V\in\mathcal{V}({\overline{y}})$ such that, for every $y\in V,$%
\begin{equation}
d(\overline{x},F^{-1}(y))\leq L\left\Vert y-{\overline{y}}\right\Vert .
\label{Lhemireg}%
\end{equation}

The infimum of $L>0$ over all the combinations $(L,V)$ for which
(\ref{Lhemireg}) holds is denoted by $\operatorname*{hemreg}F(\overline
{x},\overline{y})$ and is called the exact hemiregularity bound of $F$ at
$(\overline{x},\overline{y}).$
\end{df}

For more details about these concepts, see \cite{ArtMord2010}, \cite{DurStr4}.
The next proposition presents the corresponding equivalences between the "at
point" notions introduced before.

\begin{pr}
\label{link_at}Let $L>0,$ $F:X\rightrightarrows Y$ and $(\overline
{x},\overline{y})\in\operatorname{Gr}F.$ Then $F$ is $L-$open at
$(\overline{x},\overline{y})$ iff $F^{-1}$ is $L^{-1}-$pseudocalm at
$(\overline{y},\overline{x})$ iff $F$ is $L^{-1}-$metrically hemiregular at
$(\overline{x},\overline{y})$. Moreover, in every of the previous situations,%
\[
(\operatorname*{plop}F(\overline{x},{\overline{y}}))^{-1}%
=\operatorname*{psdclm}F^{-1}({\overline{y},}\overline{x}%
)=\operatorname*{hemreg}F(\overline{x},{\overline{y}}).
\]

\end{pr}

\bigskip

Recall that $\mathcal{L}(X,Y)$ denotes the normed vector space of linear
bounded operators acting between $X$ and $Y$. If $A\in\mathcal{L}(X,Y),$ then
the "at" and "around point" notions do coincide. In fact, $A$ is metrically
regular around every $x\in X$ iff $A$ is metrically hemiregular at every $x\in
X$ iff $A$ is open with linear rate around every $x\in X$ iff $A$ is open with
linear rate at every $x\in X$ iff $A$ is surjective. Moreover, in every of
these cases we have%
\[
\operatorname*{hemreg}A=\operatorname*{reg}A=(\operatorname*{plop}%
A)^{-1}=(\operatorname*{lop}A)^{-1}=\left\Vert (A^{\ast})^{-1}\right\Vert ,
\]

\noindent where $A^{\ast}\in\mathcal{L}(Y^{\ast},X^{\ast})$ denotes the
adjoint operator and $\operatorname*{hemreg}A,$ $\operatorname*{reg}A,$
$\operatorname*{plop}A$ and $\operatorname*{lop}A$ are common for all the
points $x\in X$ (see, for more details, \cite[Proposition 5.2]{ArtMord2010}).

\bigskip

Finally, we introduce the corresponding partial notions of linear openness,
metric regularity and Lipschitz-like property around the reference point for a
parametric set-valued map.

\begin{df}
Let $L>0,$ $F:X\times P\rightrightarrows Y$ be a multifunction, $((\overline
{x},\overline{p}),\overline{y})\in\operatorname{Gr}F$ and for every $p\in P,$
denote $F_{p}(\cdot):=F(\cdot,p).$

(i) $F$ is said to be open at linear rate $L>0,$ or $L-$open, with respect to
$x$ uniformly in $p$ around $((\overline{x},\overline{p}),\overline{y})$ if
there exist a positive number $\varepsilon>0$ and some neighborhoods
$U\in\mathcal{V}(\overline{x}),$ $V\in\mathcal{V}(\overline{p}),$
$W\in\mathcal{V}(\overline{y})$ such that, for every $\rho\in(0,\varepsilon),$
every $p\in V$ and every $(x,y)\in\operatorname*{Gr}F_{p}\cap\lbrack U\times
W],$%
\begin{equation}
B(y,\rho L)\subset F_{p}(B(x,\rho)), \label{pLopen}%
\end{equation}

The supremum of $L>0$ over all the combinations $(L,U,V,W,\varepsilon)$ for
which (\ref{pLopen}) holds is denoted by $\widehat{\operatorname*{lop}}%
_{x}F((\overline{x},\overline{p}),\overline{y})$ and is called the exact
linear openness bound, or the exact covering bound of $F$ in $x$ around
$((\overline{x},\overline{p}),\overline{y}).$

(ii) $F$ is said to be Lipschitz-like, or has Aubin property, with respect to
$x$ uniformly in $p$ around $((\overline{x},\overline{p}),\overline{y})$ with
constant $L>0$ if there exist some neighborhoods $U\in\mathcal{V}(\overline
{x}),$ $V\in\mathcal{V}(\overline{p}),$ $W\in\mathcal{V}(\overline{y})$ such
that, for every $x,u\in U$ and every $p\in V,$%
\begin{equation}
F_{p}(x)\cap W\subset F_{p}(u)+L\left\Vert x-u\right\Vert \mathbb{D}_{Y}.
\label{pLLip_like}%
\end{equation}

The infimum of $L>0$ over all the combinations $(L,U,V,W)$ for which
(\ref{pLLip_like}) holds is denoted by $\widehat{\operatorname*{lip}}%
_{x}F((\overline{x},\overline{p}),\overline{y})$ and is called the exact
Lipschitz bound of $F$ in $x$ around $((\overline{x},\overline{p}%
),\overline{y}).$

(iii) $F$ is said to be metrically regular with respect to $x$ uniformly in
$p$ around $((\overline{x},\overline{p}),\overline{y})$ with constant $L>0$ if
there exist some neighborhoods $U\in\mathcal{V}(\overline{x}),$ $V\in
\mathcal{V}(\overline{p}),$ $W\in\mathcal{V}(\overline{y})$ such that, for
every $(x,p,y)\in U\times V\times W,$%
\begin{equation}
d(x,F_{p}^{-1}(y))\leq Ld(y,F_{p}(x)). \label{pLmreg}%
\end{equation}

The infimum of $L>0$ over all the combinations $(L,U,V,W)$ for which
(\ref{pLmreg}) holds is denoted by $\widehat{\operatorname*{reg}}%
_{x}F((\overline{x},\overline{p}),\overline{y})$ and is called the exact
regularity bound of $F$ in $x$ around $((\overline{x},\overline{p}%
),\overline{y}).$
\end{df}

Similarly, one can define the notions of linear openness, metric regularity
and Lipschitz-like property with respect to $p$ uniformly in $x,$ and the
corresponding exact bounds.

\section{Linear openness of compositions}

We start the main section of the paper with a refinement of a result
previously given in \cite{DurStr4}. Here, we present some conclusions in a
slightly different form but without inner semicontinuity assumptions (see as
well the comments after the proof).

\begin{thm}
\label{impl}Let $X,P$ be metric spaces, $Y$ be a normed vector space,
$H:X\times P\rightrightarrows Y$ be a set-valued map and $(\overline
{x},\overline{p},0)\in\operatorname{Gr}H$. Denote by $H_{p}(\cdot
):=H(\cdot,p),$ $H_{x}(\cdot):=H(x,\cdot).$

(i) If $H$ is open with linear rate $c>0$ with respect to $x$ uniformly in $p$
around $(\overline{x},\overline{p},0)$, then there exist $\alpha,\beta
,\gamma>0$ such that, for every $(x,p)\in B(\overline{x},\alpha)\times
B(\overline{p},\beta),$%
\begin{equation}
d(x,S(p))\leq c^{-1}d(0,H(x,p)\cap B(0,\gamma)). \label{xSpV}%
\end{equation}

Suppose, in addition, that $Y$ is a normed vector space and $H$ is
Lipschitz-like with respect to $p$ uniformly in $x$ around $(\overline
{x},\overline{p},0).$ Then $S$ is Lipschitz-like around $(\overline
{p},\overline{x})$ and%
\begin{equation}
\operatorname*{lip}S(\overline{p},\overline{x})\leq c^{-1}\widehat
{\operatorname*{lip}}_{p}H((\overline{x},\overline{p}),0).\label{lipS}%
\end{equation}

(ii) If $H$ is open with linear rate $c>0$ with respect to $p$ uniformly in
$x$ around $(\overline{x},\overline{p},0)$, then there exist $\alpha
,\beta,\gamma>0$ such that, for every $(x,p)\in B(\overline{x},\alpha)\times
B(\overline{p},\beta),$%
\begin{equation}
d(p,S^{-1}(x))\leq c^{-1}d(0,H(x,p)\cap B(0,\gamma)). \label{pSxV}%
\end{equation}

If, moreover, $H$ is Lipschitz-like with respect to $x$ uniformly in $p$
around $(\overline{x},\overline{p},0),$ then $S$ is metrically regular around
$(\overline{p},\overline{x})$ and%
\begin{equation}
\operatorname*{reg}S(\overline{p},\overline{x})\leq c^{-1}\widehat
{\operatorname*{lip}}_{x}H((\overline{x},\overline{p}),0). \label{regS}%
\end{equation}

\end{thm}

\noindent\textbf{Proof. }We will prove only the first item, because for the
second one it suffices to observe that, defining the multifunction
$T:=S^{-1},$ the proof is completely symmetrical, using $T$ instead of $S$.
Moreover, using Proposition \ref{link_around}, we know that
$\operatorname*{reg}S(\overline{p},\overline{x})=\operatorname*{lip}%
T(\overline{x},\overline{p})$ and then (\ref{regS}) follows from (\ref{lipS}).

For the (i) item, we know that there exist $r,s,t,c,\varepsilon>0$ such that,
for every $\rho\in(0,\varepsilon),$ every $p\in B(\overline{p},t)$ and every
$(x,y)\in\operatorname{Gr}H_{p}\cap\lbrack B(\overline{x},r)\times B(0,s)],$%
\[
B(y,c\rho)\subset H_{p}(B(x,\rho)).
\]

Take now $\rho\in(0,\min\{\varepsilon,c^{-1}s\}).$ Set $\alpha:=r,$
$\beta:=t,$ $\gamma:=c\rho$ and fix arbitrary $(x,p)\in B(\overline{x}%
,\alpha)\times B(\overline{p},\beta).$ If $H(x,p)\cap B(0,c\rho)=\emptyset,$
then $d(0,H(x,p)\cap B(0,\gamma))=+\infty$ and (\ref{xSpV}) trivially holds.
Suppose next that $H(x,p)\cap B(0,c\rho)\not =\emptyset.$ If $0\in H(x,p),$
then $0\in H(x,p)\cap B(0,c\rho),$ and, again, (\ref{xSpV}) trivially holds.
Suppose now that $0\not \in H(x,p)\cap B(0,c\rho).$ Then for every $\xi>0,$
there exists $y_{\xi}\in H(x,p)\cap B(0,c\rho)$ such that%
\[
\left\Vert y_{\xi}\right\Vert <d(0,H(x,p)\cap B(0,c\rho))+\xi.
\]

Because $d(0,H(x,p)\cap B(0,c\rho))<c\rho,$ we can choose $\xi$ sufficiently
small such that $d(0,H(x,p)\cap B(0,c\rho))+\xi<c\rho.$ Consequently,%
\begin{equation}
0\in B(y_{\xi},d(0,H(x,p)\cap B(0,c\rho))+\xi)\subset B(y_{\xi},c\rho).
\label{0inB}%
\end{equation}

Observe now that $x\in B(\overline{x},r),$ $p\in B(\overline{p},t),$ $y_{\xi
}\in B(0,d(0,H(x,p)\cap B(0,c\rho))+\xi)\subset B(0,c\rho)\subset B(0,s),$
$y_{\xi}\in H(x,p)$ and denote $\rho_{0}:=c^{-1}(d(0,H(x,p)\cap B(0,c\rho
))+\xi)<\rho<\varepsilon.$

But we know that%
\[
B(y_{\xi},c\rho_{0})\subset H_{p}(B(x,\rho_{0})),
\]

\noindent hence, using also (\ref{0inB}), one obtains that there exists
$x_{0}\in B(x,\rho_{0})$ such that $0\in H(x_{0},p),$ which is equivalent to
$x_{0}\in S(p).$ Then%
\[
d(x,S(p))\leq d(x,x_{0})<\rho_{0}=c^{-1}(d(0,H(x,p)\cap B(0,c\rho))+\xi).
\]

Making $\xi\rightarrow0,$ we obtain (\ref{xSpV}).

Suppose now that $H$ is Lipschitz-like with respect to $p$ uniformly in $x$
around $(\overline{x},\overline{p},0).$ Then there exist $l,a,b,\tau>0$ such
that $\tau<c\rho$ and for every $x\in B(\overline{x},a)$ and every
$p_{1},p_{2}\in B(\overline{p},b),$%
\begin{equation}
H(x,p_{1})\cap D(0,\tau)\subset H(x,p_{2})+ld(p_{1},p_{2})\mathbb{D}_{Y}.
\label{lipHp}%
\end{equation}

Take $\overline{\alpha}:=\min\{a,\alpha\},$ $\overline{\beta}=\min
\{b,\beta,(2l)^{-1}\tau\},$ $p_{1},p_{2}\in B(\overline{p},\overline{\beta})$
and $x\in S(p_{1})\cap B(\overline{x},\overline{\alpha}).$ Then $0\in
H(x,p_{1})\cap D(0,\tau),$ whence, using (\ref{lipHp}), there exists
$y^{\prime\prime}\in\mathbb{D}_{Y}$ such that $y^{\prime}:=l\cdot
d(p_{1},p_{2})y^{\prime\prime}\in H(x,p_{2})$ with $\left\Vert y^{\prime
}\right\Vert \leq l\cdot\lbrack d(p_{1},\overline{p})+d(\overline{p}%
,p_{2})]\leq\tau<c\rho.$ Hence, $y^{\prime}\in H(x,p_{2})\cap B(0,\gamma),$ so
using (\ref{xSpV}), we get that%
\[
d(x,S(p_{2}))\leq c^{-1}d(0,H(x,p_{2})\cap B(0,\gamma))\leq c^{-1}\left\Vert
y^{\prime}\right\Vert \leq c^{-1}ld(p_{1},p_{2}).
\]

Consequently, because $l$ can be chosen arbitrarily close to $\widehat
{\operatorname*{lip}}_{p}H((\overline{x},\overline{p}),0)),$ it follows that
$S$ is Lipschitz-like around $(\overline{p},\overline{x})$ and
$\operatorname*{lip}S(\overline{p},\overline{x})\leq c^{-1}\widehat
{\operatorname*{lip}}_{p}H((\overline{x},\overline{p}),0)).$ The proof is now
complete.$\hfill\square$

\begin{rmk}
If, in addition to the assumptions of Theorem \ref{impl}, $H$ is inner
semicontinuous at $(\overline{x},\overline{p},0),$ then the relation
(\ref{xSpV}) becomes%
\[
d(x,S(p))\leq c^{-1}d(0,H(x,p)).
\]

Also, as one can see from the precedent proof, $P$ can be taken to be just
topological space (see, for more details, \cite[Theorem 3.2]{DurStr3}).
\end{rmk}

We present the next result as a by-product of Theorem \ref{impl}. Note that
different versions of the second part of the following lemma are done in
\cite[Theorem 3.5]{ArtMord2010} (with functions instead of multifunctions) and
in \cite[Lemma 2]{Ioffe2010b}. Here we also obtain some extra conclusions for
a general situation. Taking into account that this lemma will be used in the
sequel, we prefer to give all the details of its proof.

\begin{lm}
\label{gama}Let $Y,Z,W$ be normed vector spaces, $G:Y\times Z\rightrightarrows
W$ be a multifunction and $({\overline{y},}\overline{z},\overline{w})\in
Y\times Z\times W$ be such that $\overline{w}\in G(\overline{y},{\overline{z}%
}).$ Consider next the implicit multifunction $\Gamma:Z\times
W\rightrightarrows Y$ defined by%
\[
\Gamma(z,w):=\{y\in Y\mid w\in G(y,z)\}.
\]

Suppose that the following conditions are satisfied:

(i) $G$ is Lipschitz-like with respect to $z$ uniformly in $y$ around
$((\overline{y},{\overline{z}),}\overline{w})$ with constant $D\geq0;$

(ii) $G$ is open at linear rate with respect to $y$ uniformly in $z$ around
$((\overline{y},{\overline{z}),}\overline{w})$ with constant $C>0.$

Then the multifunction $H:Y\times Z\times W\rightrightarrows W,$ given by
\[
H(y,(z,w)):=G(y,z)-w\text{ for every }(y,z,w)\in Y\times Z\times W,
\]
is Lipschitz-like with respect to $(z,w)$ uniformly in $y$ around
$(\overline{y},({\overline{z},}\overline{w}),0).$

Moreover, there exists $\gamma>0$ such that for every $(z,w),(z^{\prime
},w^{\prime})\in D({\overline{z},}\gamma)\times D(\overline{w},\gamma),$ and
for every $\delta>0,$
\begin{equation}
\Gamma(z,w)\cap D(\overline{y},\gamma)\subset\Gamma(z^{\prime},w^{\prime
})+\frac{1+\delta}{C}(D\left\Vert z-z^{\prime}\right\Vert +\left\Vert
w-w^{\prime}\right\Vert )\mathbb{D}_{Y}. \label{lipGam}%
\end{equation}

\end{lm}

\noindent\textbf{Proof.} Define $P:=Z\times W$ and then $H:Y\times
P\rightrightarrows W.$ Observe also that%
\[
\Gamma(z,w)=\{y\in Y\mid0\in H(y,(z,w))\},
\]
and denote $H_{(z,w)}(\cdot):=H(\cdot,(z,w)).$ Because for every $z$ close to
${\overline{z}}$ we know from (ii) that $G_{z}$ is $C-$open at points from its
graph around $(\overline{y},\overline{w}),$ we can conclude that there exists
$\varepsilon>0$ such that for every $\rho\in(0,\varepsilon),$ every $z\in
B({\overline{z},}\varepsilon)$ and every $(y,w^{\prime})\in\operatorname{Gr}%
G_{z}\cap\lbrack B(\overline{y},\varepsilon)\times B(\overline{w}%
,\varepsilon)],$
\[
B(w^{\prime},C\rho)\subset G_{z}(B(y,\rho)).
\]

Hence, for every $(z,w)\in B({\overline{z},}\varepsilon)\times B({\overline
{w},2}^{-1}\varepsilon)$ and every $(y,u)\in\operatorname{Gr}H_{(z,w)}%
\cap\lbrack B(\overline{y},\varepsilon)\times B({0,2}^{-1}\varepsilon)],$
setting $w^{\prime}:=u+w\in B(\overline{w},\varepsilon),$ we get that%
\[
B(u,C^{-1}\rho)=B(w^{\prime},C\rho)-w\subset G_{z}(B(y,\rho))-w=H_{(z,w)}%
(B(y,\rho)).
\]
But this shows, applying Theorem \ref{impl}, that there exist $\beta>0$ such
that, for every $(y,z,w)\in B(\overline{y},\beta)\times B(\overline{z}%
,\beta)\times B(\overline{w},\beta),$%
\begin{equation}
d(y,\Gamma(z,w))\leq C^{-1}d(0,H(y,(z,w))\cap B(0,\beta)). \label{relGa}%
\end{equation}

We want to prove that there exists $c>0$ such that $(D+1)c<\beta$, for every
$y\in B(\overline{y},c),$ $(z,w),(z^{\prime},w^{\prime})\in B({\overline{z}%
,}c)\times B(\overline{w},c),$%

\begin{equation}
H(y,(z,w))\cap D(0,c)\subset H(y,(z^{\prime},w^{\prime}))+(D\left\Vert
z-z^{\prime}\right\Vert +\left\Vert w-w^{\prime}\right\Vert )\mathbb{D}_{Z}.
\label{relH}%
\end{equation}

In particular, we will prove the first part of the conclusion.

Because of (i), we know that there exists $a>0$ such that for every $y\in
B(\overline{y},a)$ and every $z,z^{\prime}\in B({\overline{z},a}),$%
\begin{equation}
G(y,z)\cap D(\overline{w},a)\subset G(y,z^{\prime})+D\left\Vert z-z^{\prime
}\right\Vert \mathbb{D}_{Z}. \label{relG}%
\end{equation}

Choose now $c\in(0,\min\{2^{-1}a,(D+1)^{-1}\beta\})$ and take arbitrary $y\in
B(\overline{y},c),$ $(z,w),(z^{\prime},w^{\prime})$ $\in B({\overline{z}%
,}c)\times B(\overline{w},c).$ Furthermore, choose $u^{\prime}\in
H(y,(z,w))\cap D(0,c).$ Then $\left\Vert u^{\prime}+w-\overline{w}\right\Vert
\leq2c<a,$ whence $u^{\prime}+w\in G(y,z)\cap D(\overline{w},a)$ and because
of (\ref{relG}), one obtains successively that
\begin{align*}
u^{\prime}+w &  \in G(y,z^{\prime})+D\left\Vert z-z^{\prime}\right\Vert
\mathbb{D}_{Z},\\
u^{\prime} &  \in G(y,z^{\prime})-w^{\prime}+(w-w^{\prime})+D\left\Vert
z-z^{\prime}\right\Vert \mathbb{D}_{Z},\\
u^{\prime} &  \in H(y,(z^{\prime},w^{\prime}))+(D\left\Vert z-z^{\prime
}\right\Vert +\left\Vert w-w^{\prime}\right\Vert )\mathbb{D}_{Z}.
\end{align*}

Take now $\gamma\in(0,\min\{\beta,2^{-1}c\}),$ $(z,w),(z^{\prime},w^{\prime
})\in B({\overline{z},}\gamma)\times B(\overline{w},\gamma)$ and $y\in
\Gamma(z,w)\cap D(\overline{y},\gamma).$ Hence, $0\in H(y,(z,w))\cap D(0,c).$
Then, using (\ref{relH}), we have that there exists $\eta\in\mathbb{D}_{Z}$
such that $(D\left\Vert z-z^{\prime}\right\Vert +\left\Vert w-w^{\prime
}\right\Vert )\eta\in H(y,(z^{\prime},w^{\prime}))\cap B(0,\beta)$ (because
$(D+1)c<\beta$). Finally, using (\ref{relGa}), we get that
\begin{align*}
d(y,\Gamma(z^{\prime},w^{\prime}))  &  \leq C^{-1}d(0,H(y,(z^{\prime
},w^{\prime}))\cap B(0,\beta))\\
&  \leq C^{-1}(D\left\Vert z-z^{\prime}\right\Vert +\left\Vert w-w^{\prime
}\right\Vert ),
\end{align*}

\noindent which completes the proof.$\hfill\square$

\bigskip

We are now in position to formulate and to prove our main result, an openness
result for a fairly general set-valued composition.

\begin{thm}
\label{op_comp}Let $X,Y,Z,W$ be Banach spaces, $F_{1}:X\rightrightarrows Y,$
$F_{2}:X\rightrightarrows Z$ and $G:Y\times Z\rightrightarrows W$ be three
multifunctions and $(\overline{x},{\overline{y}},\overline{z},\overline{w})\in
X\times Y\times Z\times W$ such that $(\overline{x},{\overline{y}}%
)\in\operatorname{Gr}F_{1},$ $(\overline{x},{\overline{z}})\in
\operatorname{Gr}F_{2}$ and $(({\overline{y},}\overline{z}),\overline{w}%
)\in\operatorname{Gr}G.$ Let $H:X\rightrightarrows W$ be given by
\[
H(x):=G(F_{1}(x),F_{2}(x))\text{ for every }x\in X
\]
and suppose that the following assumptions are satisfied:

(i) $\operatorname{Gr}F_{1},$ $\operatorname{Gr}F_{2}$ and $\operatorname{Gr}%
G$ are locally closed around $(\overline{x},{\overline{y}}),$ $(\overline
{x},\overline{z})$ and $(({\overline{y},}\overline{z}),\overline{w});$

(ii) $F_{1}$ is open at linear rate $L>0$ around $(\overline{x},{\overline{y}%
});$

(iii) $F_{2}$ is Lipschitz-like around $(\overline{x},{\overline{z}})$ with
constant $M>0;$

(iv) $G$ is open at linear rate with respect to $y$ uniformly in $z$ around
$((\overline{y},{\overline{z}),}\overline{w})$ with constant $C>0;$

(v) $G$ is Lipschitz-like with respect to $z$ uniformly in $y$ around
$((\overline{y},{\overline{z}),}\overline{w})$ with constant $D\geq0;$

(vi) $LC-MD>0.$

Then there exists $\varepsilon>0$ such that, for every $\rho\in(0,\varepsilon
),$%
\[
B(\overline{w},(LC-MD)\rho)\subset H(B(\overline{x},\rho)).
\]

Moreover, for every $\rho\in(0,2^{-1}\varepsilon)$ and every $(x,y,z,w)\in
B(\overline{x},2^{-1}\varepsilon)\times B({\overline{y},2}^{-1}\varepsilon
)\times{B(}\overline{z},2^{-1}\varepsilon)\times B(\overline{w},2^{-1}%
\varepsilon)$ such that $(y,z)\in F_{1}(x)\times F_{2}(x)$ and $w\in G(y,z),$
\[
B(w,(LC-MD)\rho)\subset H(B(x,\rho)).
\]

\end{thm}

\noindent\textbf{Proof.} Using (i), one can find $\alpha>0$ such that
$\operatorname{Gr}F_{1}\cap\operatorname*{cl}[B(\overline{x},\alpha)\times
B({\overline{y}},\alpha)],$ $\operatorname{Gr}F_{2}\cap\operatorname*{cl}%
[B(\overline{x},\alpha)\times B(\overline{z},\alpha)]$ and $\operatorname{Gr}%
G\cap\operatorname*{cl}[B(\overline{y},\alpha)\times B(\overline{z}%
,\alpha)\times B(\overline{w},\alpha)]$ are closed. Also, using (ii) and
(iii), there exist $\beta>0$ such that, for every $(x^{\prime},y^{\prime}%
)\in\operatorname{Gr}F_{1}\cap\lbrack B(\overline{x},\beta)\times
B({\overline{y}},\beta)],$ $F_{1}$ is $L-$open at $(x^{\prime},y^{\prime})$
and for every $(v^{\prime},u^{\prime})\in\operatorname{Gr}F_{2}^{-1}%
\cap\lbrack B(\overline{z},\beta)\times B(\overline{x},\beta)],$ $F_{2}^{-1}$
is $M^{-1}-$open at $(v^{\prime},u^{\prime}).$ Finally, using (iv) and (v), we
can find $\gamma>0$ such that for every $(z,w),(z^{\prime},w^{\prime})\in
B({\overline{z},}\gamma)\times B(\overline{w},\gamma)$ and every $\delta>0$%
\begin{equation}
\Gamma(z,w)\cap B(\overline{y},\gamma)\subset\Gamma(z^{\prime},w^{\prime
})+\frac{1+\delta}{C}(D\left\Vert z-z^{\prime}\right\Vert +\left\Vert
w-w^{\prime}\right\Vert )\mathbb{B}_{Y}. \label{relGam}%
\end{equation}

Without loosing the generality, using (v), we can suppose also that for every
$y\in B(\overline{y},\gamma)$ and every $z,z^{\prime}\in B({\overline{z}%
,}\gamma),$%
\begin{equation}
G(y,z)\cap B(\overline{w},\gamma)\subset G(y,z^{\prime})+D\left\Vert
z-z^{\prime}\right\Vert \mathbb{D}_{Y}. \label{lipG}%
\end{equation}

Fix $\varepsilon:=\min\{\alpha,L^{-1}\alpha,M^{-1}\alpha,(LC+MD)^{-1}%
\alpha,\beta,2^{-1}L^{-1}\beta,2^{-1}M^{-1}\beta,\gamma,2^{-1}L^{-1}%
\gamma,2^{-1}M^{-1}\gamma,$ $2^{-1}(LC+MD)^{-1}\gamma\}$ and take $\rho
\in(0,\varepsilon).$

Define the multifunction $(F_{1},F_{2}):X\rightrightarrows Y\times Z$ by
$(F_{1},F_{2})(x):=F_{1}(x)\times F_{2}(x)$ and observe that $(\overline
{x},{\overline{y}},\overline{z})\in\operatorname*{Gr}(F_{1},F_{2}).$ Because
of the choice of $\varepsilon,$ we have that the set $\Omega\cap
\operatorname*{cl}A$ is closed, where%

\begin{align}
A  &  :=B(\overline{x},\rho)\times B({\overline{y}},L\rho)\times
B(\overline{z},M\rho)\times B(\overline{w},(LC+MD)\rho),\text{ and}\label{A}\\
\Omega &  :=\{(x,y,z,w)\in X\times Y\times Z\times W\mid(y,z)\in(F_{1}%
,F_{2})(x)\text{ and }w\in G(y,z)\}. \label{omega}%
\end{align}

Take $u\in B(\overline{w},(LC-MD)\rho).$ We must prove that $u\in
H(B(\overline{x},\rho)).$ There exists $\tau\in(0,1)$ such that $\left\Vert
u-\overline{w}\right\Vert <\tau(LC-MD)\rho.$ Endow the space $X\times Y\times
Z\times W$ with the norm%
\[
\left\Vert (p,q,r,s)\right\Vert _{0}:=\tau(LC-MD)\max\{\left\Vert p\right\Vert
,L^{-1}\left\Vert q\right\Vert ,M^{-1}\left\Vert r\right\Vert ,(LC+MD)^{-1}%
\left\Vert s\right\Vert \}
\]

\noindent and apply the Ekeland variational principle to the function
$h:\Omega\cap\operatorname*{cl}A\rightarrow\mathbb{R}_{+},$%
\[
h(p,q,r,s):=\left\Vert u-s\right\Vert .
\]

Then one can find a point $(a,b,c,d)\in\Omega\cap\operatorname*{cl}A$ such
that%
\begin{equation}
\left\Vert u-d\right\Vert \leq\left\Vert u-\overline{w}\right\Vert -\left\Vert
(a,b,c,d)-(\overline{x},{\overline{y}},\overline{z},\overline{w})\right\Vert
_{0} \label{ek1}%
\end{equation}

\noindent and%
\begin{equation}
\left\Vert u-d\right\Vert \leq\left\Vert u-\overline{w}\right\Vert +\left\Vert
(a,b,c,d)-(p,q,r,s)\right\Vert _{0},\text{ }\forall(p,q,r,s)\in\Omega
\cap\operatorname*{cl}A. \label{ek2}%
\end{equation}

From (\ref{ek1}) we have that%
\begin{align*}
&  \tau(LC-MD)\max\{\left\Vert a-\overline{x}\right\Vert ,L^{-1}\left\Vert
b-{\overline{y}}\right\Vert ,M^{-1}\left\Vert c-\overline{z}\right\Vert
,(LC+MD)^{-1}\left\Vert d-\overline{w}\right\Vert \}\\
&  =\left\Vert (a,b,c,d)-(\overline{x},{\overline{y}},\overline{z}%
,\overline{w})\right\Vert _{0}\leq\left\Vert u-\overline{w}\right\Vert
<\tau(LC-MD)\rho,
\end{align*}

\noindent hence $(a,b,c,d)\in A,$ and, in particular, $a\in B(\overline
{x},\rho).$

If $u=d,$ then $u\in H(a)\subset H(B(\overline{x},\rho))$ and the desired
assertion is proved.

We want to show that $u=d$ is the sole possible situation. For this, suppose
by means of contradiction that $u\not =d.$ Fix $\omega>0$ such that
$LC-\omega>0$ and define next%
\[
v:=(LC-\omega)\left\Vert u-d\right\Vert ^{-1}(u-d).
\]

\noindent Take arbitrary $\zeta\in(0,2^{-1}\rho).$ Remark that, from the
choice of $\varepsilon,$ $b\in B({\overline{y}},\gamma),c\in B(\overline
{z},\gamma)$ and $d,d+\zeta v\in B(\overline{w},\gamma).$ Hence, using
(\ref{relGam}) for $\delta\in(0,(LC-\omega)^{-1}\omega),$ we get that
\[
b\in\Gamma(c,d)\cap B({\overline{y}},\gamma)\subset\Gamma(c,d+\zeta
v)+(L-\varepsilon^{\prime})\zeta\mathbb{D}_{Y},
\]

\noindent where $\varepsilon^{\prime}:=C^{-1}(\omega-\delta(LC-\omega
))\in(0,L).$ Therefore, there exists $q\in\mathbb{D}_{Y}$ such that
$b+(L-\varepsilon^{\prime})\zeta q\in\Gamma(c,d+\zeta v),$ or, equivalently,
$d+\zeta v\in G(b+(L-\varepsilon^{\prime})\zeta q,c).$

Now, from the $L-$ openness of $F_{1}$ at $(a,b)\in\operatorname{Gr}F_{1}%
\cap\lbrack B(\overline{x},\beta)\times B({\overline{y}},\beta)],$ there
exists $\varepsilon_{0}<2^{-1}\rho$ such that, for every $\zeta\in
(0,\varepsilon_{0})$,%
\begin{equation}
b+(L-\varepsilon^{\prime})\zeta q\in B(b,L\zeta)\subset F_{1}(B(a,\zeta)).
\label{F1}%
\end{equation}

Consequently, there exists $p$ with $\left\Vert p\right\Vert <1$ such that
$b+(L-\varepsilon^{\prime})\zeta q\in F_{1}(a+\zeta p).$

Also, using the $M^{-1}-$ openness of $F_{2}^{-1}$ at $(c,a)\in
\operatorname{Gr}F_{2}^{-1}\cap\lbrack B(\overline{z},\beta)\times
B(\overline{x},\beta)],$ we can find $\varepsilon_{1}<\varepsilon_{0}$ such
that for every $\zeta\in(0,\varepsilon_{1}),$%
\begin{equation}
B(a,\zeta)\subset F_{2}^{-1}(B(c,M\zeta)). \label{F2}%
\end{equation}

Hence, one can find $e\in B(c,M\zeta)$ such that $a+\zeta p\in F_{2}^{-1}(e)$
or, equivalently, $e\in F_{2}(a+\zeta p).$ Because we can write $e=c+M\zeta r$
with $\left\Vert r\right\Vert <1$, we finally have that $(a+\zeta
p,b+(L-\varepsilon^{\prime})\zeta q,c+M\zeta r)\in\operatorname{Gr}%
(F_{1},F_{2}).$

From the choice of $\varepsilon,$ we know that $c,c+M\zeta r\in B(\overline
{z},\gamma)$ and $b+(L-\varepsilon^{\prime})\zeta q\in B({\overline{y},\gamma
}).$ Using now (\ref{lipG}), we get that for every $\zeta\in(0,\varepsilon
_{1})$,
\begin{align*}
d+\zeta v  &  \in G(b+(L-\varepsilon^{\prime})\zeta q,c)\cap B(\overline
{w},\gamma)\\
&  \subset G(b+(L-\varepsilon^{\prime})\zeta q,c+M\zeta r)+MD\zeta
\mathbb{D}_{W},
\end{align*}

\noindent so there exist $s\in\mathbb{D}_{W}$ such that $d+\zeta v+MD\zeta
s\in G(b+(L-\varepsilon^{\prime})\zeta q,c+M\zeta r).$

In conclusion, for every $\zeta\in(0,\varepsilon_{1})$, one can find
$(p,q,r,s)\in(\mathbb{B}_{X},\mathbb{D}_{Y},\mathbb{B}_{Z},\mathbb{D}_{W})$
such that $(a+\zeta p,b+(L-\varepsilon^{\prime})\zeta q,c+M\zeta r,d+\zeta
v+MD\zeta s)\in\Omega\cap A.$We use now (\ref{ek2}) to obtain that%
\begin{align}
\left\Vert u-d\right\Vert  &  \leq\left\Vert u-(d+\zeta v+MD\zeta
s)\right\Vert +\zeta\left\Vert (p,(L-\varepsilon^{\prime}%
)q,Mr,v+MDs)\right\Vert _{0}\label{rel}\\
&  \leq\left\Vert u-d-\zeta v\right\Vert +MD\zeta+\zeta\left\Vert
(p,(L-\varepsilon^{\prime})q,Mr,v+MDs)\right\Vert _{0}.\nonumber
\end{align}

But%
\[
\left\Vert u-d-\zeta v\right\Vert =\left\vert \left\Vert u-d\right\Vert
-\zeta(LC-\omega)\right\vert .
\]

Set $\varepsilon_{2}:=\min\{\varepsilon_{1},(LC-\varepsilon)^{-1}\left\Vert
u-d\right\Vert \}.$ Then for every $\zeta\in(0,\varepsilon_{2}),$ one obtains
successively from (\ref{rel}) that%
\begin{align*}
\left\Vert u-d\right\Vert  &  \leq\left\Vert u-d\right\Vert -\zeta
(LC-\omega)+MD\zeta+\zeta\left\Vert (p,(L-\varepsilon^{\prime}%
)q,Mr,v+MDs)\right\Vert _{0},\\
LC-MD-\omega &  \leq\tau(LC-MD)\max\{\left\Vert p\right\Vert ,L^{-1}\left\Vert
(L-\varepsilon^{\prime})q\right\Vert ,M^{-1}\left\Vert Mr\right\Vert
,(LC+MD)^{-1}\left\Vert v+MDs\right\Vert \},\\
LC-MD-\omega &  \leq\tau(LC-MD).
\end{align*}

\noindent Passing to the limit when $\omega\rightarrow0,$ we get that
$1\leq\tau,$ which is a contradiction. The proof of the first part is now complete.

For the second part, the proof is similar. Suppose that the constants
$\alpha,\beta,\gamma,\varepsilon>0$ are chosen as above, take $\rho
\in(0,2^{-1}\varepsilon)$ and $(x,y,z,w)\in B(\overline{x},2^{-1}%
\varepsilon)\times B({\overline{y},2}^{-1}\varepsilon)\times{B(}\overline
{z},2^{-1}\varepsilon)\times B(\overline{w},2^{-1}\omega)$ such that $(y,z)\in
F_{1}(x)\times F_{2}(x)$ and $w\in G(y,z),$ and define $A_{1}:=B(x,\rho)\times
B({y},L\rho)\times B(z,M\rho)\times B(w,(LC+MD)\rho).$ Then $\Omega
\cap\operatorname*{cl}A_{1}$ is again closed, because $A_{1}\subset
B(\overline{x},\omega)\times B({\overline{y}},\omega)\times B(\overline
{z},\omega)\times B(\overline{w},\omega)$ and $\omega<\alpha.$ The rest of the
proof is the same as before, observing that because $(a,b,c,d)\in A_{1},$
their small perturbations for $\zeta>0$ sufficiently small remain in the
desired balls, such that one can apply relations (\ref{relGam}), (\ref{lipG}),
(\ref{F1}), (\ref{F2}).$\hfill\square$

\bigskip

As one can see, the first conclusion in Theorem \ref{op_comp} is an "at point"
openness property for $H$, while the second conclusion is not a genuine
"around point" openness property, because we cannot guarantee that the
property holds for any point $(x,w)\in\operatorname*{Gr}H$ close to
$(\overline{x},\overline{w}).$ We were able to prove that $B(w,(LC-MD)\rho
)\subset H(B(x,\rho))$ exactly because the "intermediate" points $y,z$ are
also close to $\overline{y}$ and $\overline{z},$ respectively. In general, it
is not possibly to take $y,z$ with these properties and this is the same
difficulty for getting local openness results for sum of multifunctions (see
also Corollary \ref{main_const} and the corresponding comments in
\cite{DonFra2010}).

Next results, obtained on the basis of Theorem \ref{op_comp}, improves and
extends the conclusions of \cite[Theorem 5]{Ioffe2010b}.

\begin{thm}
\label{op_comp_part}Let $X,Y,Z$ be Banach spaces, $F:X\rightrightarrows Y$ and
$G:X\times Y\rightrightarrows Z$ be two multifunctions, and $(\overline
{x},{\overline{y}},\overline{z})\in X\times Y\times Z$ such that
$(\overline{x},{\overline{y}})\in\operatorname{Gr}F,$ $(({\overline{x}%
,}\overline{y}),\overline{z})\in\operatorname{Gr}G.$ Suppose that
$\operatorname{Gr}F$ and $\operatorname{Gr}G$ are locally closed around
$(\overline{x},{\overline{y}})$ and $(({\overline{x},}\overline{y}%
),\overline{z}),$ respectively. Let $\Phi:X\rightrightarrows Z$ be given by%
\[
\Phi(x):=G(x,F(x))\text{ for every }x\in X.
\]

(A) Moreover, suppose that:

(i) $G$ is Lipschitz-like with respect to $x$ uniformly in $y$ around
$((\overline{x},{\overline{y}),}\overline{z})$ with constant $D\geq0;$

(ii) $G$ is open at linear rate with respect to $y$ uniformly in $x$ around
$((\overline{x},{\overline{y}),}\overline{z})$ with constant $C>0;$

(iii) $F$ is open at linear rate $L>0$ around $(\overline{x},{\overline{y}});$

(iv) $LC-D>0.$

Then there exists $\varepsilon>0$ such that, for every $\rho\in(0,\varepsilon
)$ and every $(x,y,z)\in B(\overline{x},\varepsilon)\times B({\overline{y}%
,}\varepsilon)\times{B(}\overline{z},\varepsilon)$ such that $y\in F(x)$ and
$z\in G(x,y),$
\begin{equation}
B(z,(LC-D)\rho)\subset\Phi(B(x,\rho)). \label{LCD}%
\end{equation}

Moreover, if $G$ satisfies the condition that for every $x\in B(\overline
{x},\varepsilon),$
\begin{equation}
G(x,y)\cap G(x,y^{\prime})=\emptyset\text{ if }y\not =y^{\prime}, \label{cond}%
\end{equation}
then the conclusion reads as follows: there exists $\varepsilon^{\prime}>0$
such that, for every $\rho\in(0,\varepsilon^{\prime})$ and every
$(x,z)\in\operatorname{Gr}\Phi\cap\lbrack B(\overline{x},\varepsilon^{\prime
})\times{B(}\overline{z},\varepsilon^{\prime})],$ the relation (\ref{LCD}) holds.

(B) Moreover, suppose that:

(i) $G$ is open at linear rate with respect to $x$ uniformly in $y$ around
$((\overline{x},{\overline{y}),}\overline{z})$ with constant $C>0;$

(ii) $G$ is Lipschitz-like with respect to $y$ uniformly in $x$ around
$((\overline{x},{\overline{y}),}\overline{z})$ with constant $D\geq0;$

(iii) $F$ is Lipschitz-like around $(\overline{x},{\overline{y}})$ with
constant $M>0;$

(iv) $C-MD>0.$

Then there exists $\varepsilon>0$ such that, for every $\rho\in(0,\varepsilon
)$ and every $(x,y,z)\in B(\overline{x},\varepsilon)\times B({\overline{y}%
,}\varepsilon)\times{B(}\overline{z},\varepsilon)$ such that $y\in F(x)$ and
$z\in G(x,y),$
\begin{equation}
B(z,(C-MD)\rho)\subset\Phi(B(x,\rho)). \label{CMD}%
\end{equation}

If, moreover, $F$ is Lipschitz around $(\overline{x},{\overline{y}})$ with
constant $M>0$ and $F(\overline{x})=\{\overline{y}\},$ then the conclusion
reads as follows: there exists $\varepsilon>0$ such that, for every $\rho
\in(0,\varepsilon)$ and every $(x,z)\in\operatorname{Gr}\Phi\cap\lbrack
B(\overline{x},\varepsilon)\times{B(}\overline{z},\varepsilon)],$ the relation
(\ref{CMD}) holds.
\end{thm}

\noindent\textbf{Proof.} The first parts of (A) and (B) easily follow from the
last conclusion of Theorem \ref{op_comp}, taking $F$ instead of $F_{1}$ or
$F_{2},$ respectively, and setting the other multifunction as the identity mapping.

For the second part of (A), consider $\Gamma:X\times Z\rightrightarrows Y$
defined by $\Gamma(x,z):=\{y\in Y\mid z\in G(x,y)\}.$ We can apply Proposition
\ref{gama} to obtain $\gamma>0$ such that for every $(x,z),(x^{\prime
},z^{\prime})\in D({\overline{x},}\gamma)\times D(\overline{z},\gamma),$ and
every $\delta>0,$%
\begin{equation}
\Gamma(x^{\prime},z^{\prime})\cap D(\overline{y},\gamma)\subset\Gamma
(x,z)+\frac{1+\delta}{C}(D\left\Vert x-x^{\prime}\right\Vert +\left\Vert
z-z^{\prime}\right\Vert )\mathbb{D}_{Y}. \label{Gamxz}%
\end{equation}

Take $\delta>0$ arbitrary small, set $\varepsilon^{\prime}\in(0,\min
\{\gamma,\varepsilon,C(D+1)^{-1}(1+\delta)^{-1}\varepsilon\}),$ and pick
$(x,z)\in\operatorname{Gr}\Phi\cap\lbrack B(\overline{x},\varepsilon^{\prime
})\times{B(}\overline{z},\varepsilon^{\prime})].$ Hence, there exists $y\in
F(x)$ such that $z\in G(x,y).$ Apply next (\ref{Gamxz}) for $(\overline{x}%
{,}\overline{z})$ instead of $(x^{\prime},z^{\prime})$ to obtain that there
exist $y^{\prime}\in\Gamma(x,z)$ such that $\left\Vert y^{\prime}%
-{\overline{y}}\right\Vert <C^{-1}(1+\delta)(D\left\Vert x-\overline
{x}\right\Vert +\left\Vert z-\overline{z}\right\Vert )<C^{-1}(1+\delta
)(D+1)\varepsilon^{\prime}<\varepsilon.$ Hence, $z\in G(x,y^{\prime})\cap
G(x,y),$ so $y=y^{\prime}\in B({\overline{y},}\varepsilon).$ The conclusion
follows now from the first parts of (A).

For the second part of (B), set $\varepsilon^{\prime}\in(0,\min\{\varepsilon
,M^{-1}\varepsilon\}),$ and pick $(x,z)\in\operatorname{Gr}\Phi\cap\lbrack
B(\overline{x},\varepsilon^{\prime})\times{B(}\overline{z},\varepsilon
^{\prime})].$ Again, there exists $y\in F(x)$ such that $z\in G(x,y).$ Using
now the additional assumptions,%
\[
y\in\{\overline{y}\}+M\left\Vert x-\overline{x}\right\Vert \mathbb{D}_{Y},
\]

\noindent whence $\left\Vert y-\overline{y}\right\Vert \leq\varepsilon,$ and
the first part applies in order to get the conclusion.$\hfill\square$

\bigskip

\begin{rmk}
The final part of the (i) item is \cite[Theorem 5]{Ioffe2010b}. Note that the
requirement of the fulfillment of (\ref{cond}) appears in \cite{Ioffe2010b}.
We observe here that it is enough to have $x$ close to ${\overline{x}}$ and,
moreover, $y^{\prime}$ could be kept close to ${\overline{y}}$.
\end{rmk}

\section{Filling the gap: metric regularity and fixed points}

In this section we point out the relations between fixed point results and
openness results. The direction that has been investigated up to now, starting
with the Arutyunov's work \cite{Arut2007}, and continued with
Dontchev-Frankowska \cite{DonFra2010} and Ioffe's \cite{Ioffe2010b} papers, is
the possibility to get openness results from fixed point assertions. Here, we
aim to emphasize that the opposite direction is also possible.

For the sake of completeness, we recall here the global version of the
Lyusternik-Graves Theorem.

\begin{thm}
Let $F:X\rightrightarrows Y$ and $G:Y\rightrightarrows X$ be two
multifunctions such that $\operatorname*{Gr}F$ and $\operatorname*{Gr}G$ are
locally closed. Suppose that $\operatorname*{Dom}(F-G^{-1})$ and
$\operatorname*{Dom}(G-F^{-1})$\ are nonempty and let $L>0$ and $M>0$ be such
that $LM>1.$ If $F$ is $L-$open at every point of its graph, and $G$ is
$M-$open at every point of its graph, then $F-G^{-1}$ is $(L-M^{-1})-$open at
every point of its graph and $G-F^{-1}$ is $(M-L^{-1})-$open at every point of
its graph.
\end{thm}

This global and fully multivalued version of the Lyusternik-Graves Theorem
appeared, for the first time, in \cite[Theorem 1]{Urs1996}. Meantime,
different variants of this result were stated in \cite{Ioffe2000},
\cite{Dmi2005}, \cite{DmiKru2008}, and the more general assumptions are that
$Y$ is a linear metric space with shift invariant metric, $X$ is metric space,
and the graphs of the involved multifunctions are complete. Another way to
prove this theorem uses an iterative procedure of Lyusternik type. As we said
before, was emphasized recently that it can also be obtained as a consequence
of fixed point theorems, as in \cite{Arut2007}, \cite{DonFra2010},
\cite{Ioffe2010b}.

\bigskip

We apply now Theorem \ref{op_comp}, for the special case where $Y=Z=W,$
$G(y,z):=y-z.$ Note that the following result was obtained in \cite{DurStr4}
as a consequence of the proof given there for global Lyusternik-Graves Theorem.

\begin{cor}
\label{main_const}Let $F_{1}:X\rightrightarrows Y$ and $F_{2}%
:X\rightrightarrows Y$ be two multifunctions and $(\overline{x},\overline
{y}_{1},\overline{y}_{2})\in X\times Y\times Y$ such that $(\overline
{x},\overline{y}_{1})\in\operatorname{Gr}F_{1}$ and $(\overline{x}%
,\overline{y}_{2})\in\operatorname{Gr}F_{2}.$ Suppose that the following
assumptions are satisfied:

(i) $\operatorname{Gr}F_{1}$ and is $\operatorname{Gr}F_{2}$ are locally
closed around $(\overline{x},\overline{y}_{1})$ and $(\overline{x}%
,\overline{y}_{2}),$ respectively;

(ii) $F_{1}$ is $l-$metrically regular around $(\overline{x},\overline{y}%
_{1});$

(iv) $F_{2}$ is $m-$Lipschitz-like around $(\overline{x},\overline{y}_{2});$

(v) $lm<1.$

Then there exists $\varepsilon>0$ such that, for every $\rho\in(0,\varepsilon
),$ $(x,y,z)\in\operatorname{Gr}(F_{1},F_{2})\cap\lbrack B(\overline
{x},\varepsilon)\times B({y}_{1}{,}\varepsilon)\times{B(}y_{2},\varepsilon)],$%
\[
B(y-z,(l^{-1}-m)\rho)\subset(F_{1}-F_{2})(B(x,\rho)).
\]

\end{cor}

\bigskip

\begin{rmk}
Let us observe that Theorem \ref{op_comp} opens several possibilities
concerning the choices of the involved multifunctions. For instance, one can
take $G(y,z):=$ $T_{1}(y)+T_{2}(z),$ where $T_{1}\in\mathcal{L}(Y,W)$ and
$T_{2}\in\mathcal{L}(Z,W)$ with $T_{1}$ surjective, and suppose that
$L\left\Vert (T_{1}^{\ast})^{-1}\right\Vert -M\left\Vert T_{2}\right\Vert
>0.\ $Then the conclusion of Theorem \ref{op_comp} applies in this special case.
\end{rmk}

\bigskip

After these comments we have the ingredients needed for closing the circle,
showing that one has the following chain of implications:

\begin{itemize}
\item Theorem $\ref{op_comp}\Rightarrow$ Corollary \ref{main_const}
$\Rightarrow$ Dontchev-Frankowska Fixed Point Theorem $\Rightarrow$ Arutyunov
Fixed Point Theorem $\Rightarrow$ global Lyusternik-Graves Theorem

\item (Proof of) global Lyusternik-Graves Theorem $\Rightarrow$ Corollary
\ref{main_const}.
\end{itemize}

\bigskip

Of course, on the first item, the first implication is given above, while the
third one is part of \cite{DonFra2010}. As we have already said, it is shown
in \cite{DurStr4}, that the proof of global Lyusternik-Graves Theorem based on
the use of Ekeland variational principle can be used to obtain Corollary
\ref{main_const}. The only implication to be proved is the second one on the
first item. This is done next.

Consider the case $\overline{y}_{1}=\overline{y}_{2}:=\overline{y}.$ Then
$0\in(F_{1}-F_{2})(\overline{x}).$ Denote
\begin{align}
S  &  :=(F_{1}-F_{2})^{-1}(0)=\{x\in X\mid0\in(F_{1}-F_{2})(x)\}\label{S}\\
&  =\{x\in X\mid F_{1}(x)\cap F_{2}(x)\neq\emptyset\}=\operatorname*{Fix}%
(F_{1}^{-1}F_{2}).\nonumber
\end{align}

\noindent Or, $S$ could be seen as the implicit multifunction associated to
$H:X\times P\rightrightarrows Y,$ $H(x,p):=F_{1}(x)-F_{2}(x)$ for every
$(x,p)\in X\times P,$ and in this case we identify the constant multifunction
$S$ with the set $S$ given by (\ref{S}).

\begin{thm}
\label{fixp}Under the assumptions of Corollary \ref{main_const}, there exist
$\alpha,\beta>0$ such that for any $x\in B(\overline{x},\alpha)$ one has that
\begin{equation}
d(x,S)\leq(l^{-1}-m)^{-1}d(F_{1}(x)\cap B(\overline{y},\beta),F_{2}(x)).
\label{diffix}%
\end{equation}

\end{thm}

\noindent\textbf{Proof. }Under the notations of Theorem \ref{main_const}, take
$\alpha,\beta>0$ such that $\alpha<m,\alpha<\varepsilon,$ $m\alpha<\beta,$
$3\beta<\varepsilon,$ $2(l^{-1}-m)^{-1}\beta<\varepsilon.$ Suppose as well
that $\alpha$ is smaller than the radius of the neighborhood of $\overline{x}$
involved in the Lipschitz-like property of $F_{2}$ at $(\overline{x}%
,\overline{y}).$ Denote $V:=B(\overline{y},\beta)$ and fix $x\in
B(\overline{x},\alpha).$ If $F_{1}(x)\cap B(\overline{y},\beta)=\emptyset$ or
$0\in F_{1}(x)\cap B(\overline{y},\beta)-F_{2}(x)$ the relation (\ref{diffix})
automatically holds. On the contrary case, for every $\mu>0$ there exists
$u_{\mu}\in F_{1}(x)\cap B(\overline{y},\beta)-F_{2}(x)$ s.t.%
\[
\left\Vert u_{\mu}\right\Vert \leq d(0,F_{1}(x)\cap B(\overline{y}%
,\beta)-F_{2}(x))+\mu.
\]
Since $F_{2}$ is Lipschitz-like around $(\overline{x},\overline{y})$ and $x\in
B(\overline{x},\alpha)$%
\[
\overline{y}\in F_{2}(\overline{x})\cap V\subset F_{2}(x)+m\left\Vert
x-\overline{x}\right\Vert \mathbb{D}_{Y},
\]
whence there is $u^{\prime}\in F_{2}(x)$ with%
\[
\left\Vert u^{\prime}-\overline{y}\right\Vert \leq m\left\Vert x-\overline
{x}\right\Vert \leq m\alpha<\beta.
\]
Therefore, $F_{2}(x)\cap V\neq\emptyset,$ whence%
\begin{align*}
\left\Vert u_{\mu}\right\Vert  &  \leq d(0,F_{1}(x)\cap B(\overline{y}%
,\beta)-F_{2}(x))+\mu\\
&  \leq d(0,F_{1}(x)\cap B(\overline{y},\beta)-F_{2}(x)\cap B(\overline
{y},\beta))+\mu\\
&  \leq2\beta+\mu.
\end{align*}
On the other hand, there are $u_{\mu}^{1}\in F_{1}(x)\cap B(\overline{y}%
,\beta)$ and $u_{\mu}^{2}\in F_{2}(x)$ with $u_{\mu}=u_{\mu}^{1}-u_{\mu}^{2}.$
Accordingly,%
\[
\left\Vert u_{\mu}^{2}-\overline{y}\right\Vert \leq\left\Vert u_{\mu}%
^{1}-\overline{y}\right\Vert +\left\Vert u_{\mu}\right\Vert \leq\beta
+2\beta+\mu=3\beta+\mu.
\]
So, if we take $\mu$ small enough one has $\left\Vert u_{\mu}^{2}-\overline
{y}\right\Vert \leq\varepsilon.$

Take $\rho_{0}=(l^{-1}-m)^{-1}d(0,F_{1}(x)\cap B(\overline{y},\beta
)-F_{2}(x))+\mu<2(l^{-1}-m)^{-1}\beta+\mu$ and hence, again for $\mu$ small
enough, $\rho_{0}<\varepsilon.$ Since $(x,u_{\mu}^{1},u_{\mu}^{2}%
)\in\operatorname{Gr}(F_{1},F_{2})\cap\lbrack B(\overline{x},\varepsilon
)\times B(\overline{y}{,}\varepsilon)\times{B(}\overline{y},\varepsilon)],$
\[
B(u_{\mu}^{1}-u_{\mu}^{2},(l^{-1}-m)\rho_{0})\subset(F_{1}-F_{2})(B(x,\rho
_{0})).
\]
Since $\left\Vert u_{\mu}\right\Vert \leq d(0,F_{1}(x)\cap B(\overline
{y},\beta)-F_{2}(x))+\mu,$ one gets%
\[
0\in B(u_{\mu}^{1}-u_{\mu}^{2},(l^{-1}-m)\rho_{0})
\]
whence there exists $x^{\prime}\in B(x,\rho_{0})$ with $0\in(F_{1}%
-F_{2})(x^{\prime}),$ whence%
\[
d(x,S)\leq(l^{-1}-m)^{-1}d(0,F_{1}(x)\cap B(\overline{y},\beta)-F_{2}(x))+\mu.
\]
Taking $\mu\rightarrow0$ one obtains the solution.$\hfill\square$

\section{Further remarks}

Theorem \ref{fixp} is the main result in \cite{DonFra2010}, and is given there
on metric spaces. Here we have obtained it on the base of Corollary
\ref{main_const}, in order to close the circle between openness of
compositions and fixed point theorems. Also, remark that a slightly different
variant of this theorem can be stated as a consequence of Theorem \ref{impl}:

\begin{pr}
Under the assumptions of Corollary \ref{main_const}, there exist $\alpha
,\beta>0$ such that for any $x\in B(\overline{x},\alpha)$ one has that
\[
d(x,S)\leq(l^{-1}-m)^{-1}d(0,(F_{1}-F_{2})(x)\cap B(0,\beta)).
\]

\end{pr}

Furthermore, the proof of this result, as the proof of the preceding theorem,
does not use the full power of Theorem \ref{impl}, which gives conclusions for
the implicit multifunction $S$ depending on the parameter $p.$ Supposing that
$F_{1}:X\times P\rightrightarrows Y,$ $F_{2}:X\rightrightarrows Y,$ and taking
$H:X\times P\rightrightarrows Y$ as $H(x,p):=F_{1}(x,p)-F_{2}(x),$ then the
implicit multifunction $S:P\rightrightarrows X$ is given by
\begin{align*}
S(p) &  :=(F_{1}(\cdot,p)-F_{2})^{-1}(0)=\{x\in X\mid0\in F_{1}(x,p)-F_{2}%
(x)\}\\
&  =\{x\in X\mid F_{1}(x,p)\cap F_{2}(x)\neq\emptyset\}=\operatorname*{Fix}%
(F_{1}(\cdot,p)^{-1}F_{2}).
\end{align*}

In this way, with a very similar proof to the one of Theorem \ref{fixp}, one
can also obtain \cite[Theorem 7]{DonFra2010}.

Of course, one can use directly Theorem \ref{op_comp} to obtain generalized
relations with $G(F_{1},F_{2})$ instead of $F_{1}-F_{2}$. Finally, observe
that almost all the openness results given before can also be stated in
parametric forms.

\end{document}